# About the trajectory synthesis to go back to nominal mode for a class of hybrid systems

**MANON Philippe, VALENTIN-ROUBINET Claire**

*Laboratoire d'Automatique et de GEnie des Procédés, UCB Lyon 1, ESCPE, UPRES-A CNRS Q5007*
*bat. 308G, 43 bd du 11 novembre 1918, 69622 VILLEURBANNE Cedex, FRANCE*
*{manon, claire}@lagep.univ-lyon1.fr*

*ABSTRACT. This paper presents basis parts of a new method to synthesize a return trajectory for a reactive process from a default mode to one of the nominal modes. The process is modeled with a hybrid automata. The purpose consists of doing a backward reachability analysis from the final state to the initial state, in the state-space. This method is applied to a batch system.*

*RESUME. Cet article présente une partie des bases théoriques d'une nouvelle méthode qui permet de synthétiser une trajectoire de retour du système réactif à un des modes nominaux depuis un mode défaillant. Le système est modélisé par un automate hybride. L'étude consiste à faire une recherche d'atteignabilité par inférence arrière dans l'espace d'état. Cette méthode est illustrée sur un système batch.*

*MOTS CLES : automate hybride, atteignabilité, champs de vecteurs linéaires, espace d'état*

# 1. Introduction

Reactive processes with equally continuous, discrete or hybrid features, may operate in different modes. They are in a nominal mode if all the process parts behave in a normal way. They are in a default mode if at least one of the different parts of the system does not behave in a normal way any more. For instance this is the case if a machine in a manufacturing process breaks down or if the level in a batch process tank oversteps the maximum or minimum allowed level. In this case, the security or the productivity decreases but the global process is not necessarily stopped. It is then a priority to find how to drive the system from its current default mode to one of its nominal mode as soon as possible.

This paper presents parts, dealing with region analytical description, of a new method which synthesizes such a trajectory for hybrid dynamical systems. The studied systems class is modeled by a hybrid automata with linear vector fields. The analysis is led in the state-space and is based on regions computations. The automata is composed of different phases. The initial phase models the system default mode and the final phase models the nominal mode the system must be driven to. The problem can be divided into two steps: a reachability analysis and a controller design.

The method presented in this paper extends Tittus's work [TIT 93], [TIT 98] on integrators hybrid systems to first order switching hybrid systems with linear vector fields. This evolution is important because we now can deal with a larger class of systems as, for example, the concentration dynamics in a batch system reactor.

Firstly, a brief example with a system in a default mode which must be driven back to its nominal mode is presented in order to clearly fix the position of the problem.

Then, the useful theory basis which is necessary to solve the problem will be developed. In section 5.2., we will explain how the sets we use can be formally described by a inequations system. The corresponding algorithm will be then briefly described. Other points of the method are presented in [MAN 99a] as the extension criteria which are necessary to solve the reachability problem.

In the last section, we will illustrate the algorithm on the example presented in section 2 and show the results.

# 2. Illustration of the problem

The method explained in this paper will be illustrated on the two tanks batch process presented figure 1. The global constraints define the minimal and maximal levels (0m and 4m) during every operating phase for the two tanks. The nominal mode limits are such that the liquid height in each tank stays between a low level (2m) and a high level (2.1m). If this constraint is violated, the system will enter a default mode. In the example, we may assume that a breakdown on the valve V1 has lead the system in a default mode. The level in each tank is then less than 0.1m, i.e. they are almost empty. The question we must solve is to decide whether or not the nominal mode is reachable from this current default mode and what is a possible control sequence which achieves it, in case it is feasible.

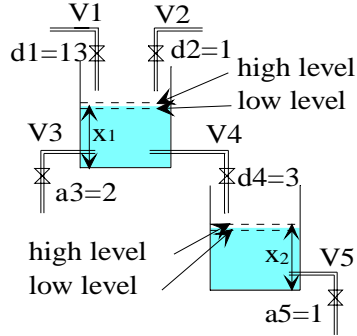

**Figure 1.** *The batch process*

## 3. The hybrid modeling

The hybrid process is modeled by a hybrid automata [TIT 98], [ASA 95]. As well as Branicky unified model [BRA 94], autonomous and controlled jumps (continuous variable discontinuity) and autonomous and controlled switching (vector field discontinuity) may be represented.

*Definition 1:*

A hybrid automata *HA* is defined by a tuple HA=($V_D$, Q, $\Sigma c$, $\mu_1$, $\mu_2$, $\mu_3$, $Q_m$):

- $V_D$ is the real valued state variables set: $V_D = \{x_i, i \in \{1, ..., n\}\}$. *x* is the continuous state vector of the system: $x=[x_1 \ldots x_n]^T \in V_D^n$. The continuous state *x* is continuously observable.
- *Q* is the possible phases set (or locations or discrete states) of the hybrid process: $Q = \{l_j, j \in \{1, ..., m\}\}$. The system global state is given by the couple $(x,l) \in V_D^n \times Q$.
- $\Sigma$ is the events set. $\Sigma = \{\delta_{i,j}, i \in \{1, ..., m\}, j \in \{1, ..., m\}/i \neq j\}$.
- $\mu_1$ is the set of *m* vector fields associated with each phase $l_j$. It models the dynamics of the state vector *x*: $\mu_1 = \{\mu_1(l_j), j \in \{1, ..., m\}\}$, that is: $\forall l_j \in Q, \mu_1(l_j)(x) = \dfrac{dx}{dt} = f(x)$ with f Lipschitz continuous.

$\mu_2$ is the set of constraints for each phase: $\mu_2 = \{\mu_2(l_i)/i \in \{1,...,m\}\}$. $\mu_2(l_j)$ is the constraint for phase $l_j$. It is a set of *l* boundaries described by the linear inequalities: $C x \leq d$ with dim(C)=lxn and dim(d)=lx1. It is not necessarily a state-space partition. A constraint $\mu_2(l_j)$ can not be violated, which means that the phase $l_j$ is left as soon as a constraint $\mu_2(l_j)$ is going to be violated. A discrete transition can be either autonomous or controlled. If the constraint $\mu_2(l_j)$ is going to be violated, a uncontrolled event occurs and the transition is autonomous, else a controlled event occurs and the transition is controlled.
- $\mu_3$ is a set of functions associated to each phase transition: $\mu_3 = \{\mu_3(l_j, \delta_{jk}, l_k) / j \in \{1,...,m\}, k \in \{1,...,m\}, j \neq k\}$. It defines real valued variables jumps when a transition occurs. $\forall (l_j, l_k) \in Q \times Q$:

　　- if the phase transition exists, then $\delta_{jk}$ is defined and
　　　　$\mu_3(l_j, \delta_{jk}, l_k) : \mathbb{R} \to \mathbb{R} / x(t^+) = \mu_3(l_j, \delta_{jk}, l_k)(x(t^-))$
　　- if the phase transition does not exist, then $\delta_{jk}$ is not defined and
　　　　$\mu_3(l_j, \delta_{jk}, l_k)$ is not defined.

　　If $\mu_3(l_j, \delta_{jk}, l_k)(x(t^-)) \neq x(t^-)$, a state jump (or state discontinuity) occurs.
- $Q_m$ is a set of marked phases, that is initial or final.　　　　　　　　　　　　Y

In order to simplify the hybrid automata graphical representation:
- $\mu_1(l_j)$ is not represented inside a phase if the vector field $\mu_1(l_j)(x)=dx/dt=0$, that is if the state-variables are constant;
- $\mu_2(l_j)$ is not represented inside a phase if it is equal to the global constraint *GC*,
- marked phases are represented by double circles.

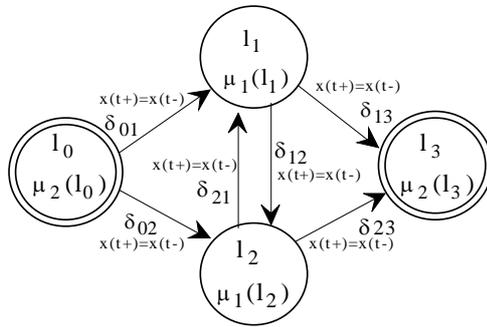

$\mu_1(l_1): \begin{cases} \dot{x}_1 = d_1 - d_4 - a_3.x_1 = 10 - 2x_1 \\ \dot{x}_2 = d_4 - a_5.x_2 = 3 - x_2 \end{cases}$

$\mu_1(l_2): \begin{cases} \dot{x}_1 = d_2 - d_4 - a_3.x_1 = -2 - 2x_1 \\ \dot{x}_2 = d_4 - a_5.x_2 = 3 - x_2 \end{cases}$

$x = \begin{bmatrix} x_1 & x_2 \end{bmatrix}^T$

$\mu_2(l_0): 0 \leq x \leq 0.1$

$\mu_2(l_3): 2 \leq x \leq 2.1$

$GC: 0 \leq x \leq 4$

**Figure 2.** *Batch process model*

All the possible behaviors of the batch process presented figure 1 are modeled figure 2. The state variables $x_1$ and $x_2$ are the liquid heights in the two tanks.

There are two types of valves. *V1, V2* and *V4* have a constant flow which is equal to $d_1$, $d_2$ and $d_4$. Valves *V3* and *V5* have a flow which is considered proportional with the level in the upstream tank ($a_3.x_1$ and $a_5.x_2$).

The operating constraints force either the valves V1, V3, V4 and V5 to be simultaneously open, which defines the *phase 1*, or the valves V2, V3, V4 and V5 to be simultaneously open, which defines the *phase 2*.

The maximum height in the tanks is 4m. The global constraint GC is then defined as: $0 \leq x \leq 4$. Every event of $\Sigma$ is controllable.

## 4. Controllability for hybrid dynamic systems (HDS): theoretical basis

If it is possible to drive the process from the initial state to the final state the system is controllable [TIT 93], [TIT 98]. Let remember the fundamental definitions for the hybrid systems controllability. They will be useful to fix the reachability algorithm.

*Definition 2:*
A hybrid system (of any dimension) is controllable if and only if there exists at least one sequence of discrete control inputs (i.e. of controlled events) that gives an acceptable trajectory along some path in the hybrid automata from any initial state to any final state in a finite time and which respects all the state constraints. Y

One or more state of the hybrid automata may not belong to the path, as well as one or more states may appear several times in the path. In this method, the reachability analysis is realized backwards, from the final state to the initial state.

The problem is to prove that a transition to any phase $l_{i+1}$ of the path does not lead to constraint violation, that means that $\mu_2(l_{i+1})$ is respected. The set of acceptable values for the continuous variables which respect $\mu_2(l_i)$ and are included in $\mu_2(l_{i+1})$ after the jump defined by the function $\mu_3(l_i, \delta_{i,i+1}, l_{i+1})$ must be calculated. This set will be named jump-set from phase $l_i$ to phase $l_{i+1}$ and will be noted $\lambda_{li \rightarrow li+1}$. Moreover, this jump-set from phase $l_i$ to phase $l_{i+1}$ may be reached by the process continuous state with the continuous dynamic evolution $\mu_1(l_i)$. Then, the jump-set can be extended to define the extended-jump-set from phase $l_i$ to phase $l_{i+1}$. It is noted $\lambda ext_{li \rightarrow li+1}$.

These concepts are defined in definition 3 and 4 and illustrated figure 3.

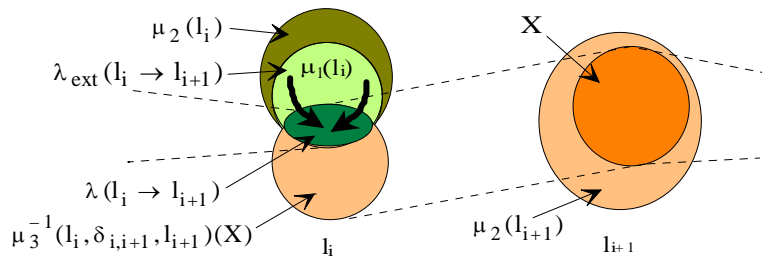

**Figure 3.** *Jump-set $\lambda_{li \rightarrow li+1}$ and extended jump-set $\lambda ext_{li \rightarrow li+1}$*

*Definition 3:*

if $l_i$ and $l_{i+1}$ are two consecutive phases of the hybrid automata, the jump-set from phase $l_i$ to phase $l_{i+1}$ ($\lambda_{li \to li+1}$) is the intersection of the phase constraints $\mu_2(l_i)$ with the inverse image region of the desired reachable set X by the inverse function $(\mu_3(l_i, \delta_{(i,i+1)}, l_{i+1}))^{-1}(X)$:

$$\forall x \in \lambda_{li \to li+1}, \mu_3(l_i, \delta_{(i,i+1)}, l_{i+1})(x) \in X$$
$$\forall x \notin \lambda_{li \to li+1}, \mu_3(l_i, \delta_{(i,i+1)}, l_{i+1})(x) \notin X$$
$$\text{i.e. } \lambda(l_i \to l_{i+1}) = \mu_2(l_i) \cap \mu_3^{-1}(l_i, \delta_{i,i+1}, l_{i+1})(X) \qquad Y$$

More detail will be given about X after the definition of the extended jump-set.

*Definition 4:*

if $l_i$ and $l_{i+1}$ are two consecutive phases of the hybrid automata, the extended-jump-set from phase $l_i$ to phase $l_{i+1}$ ($\lambda ext_{li \to li+1}$), is the region inside the phase constraint $\mu_2(l_i)$ such that $\lambda_{li \to li+1}$ is reachable by respecting the dynamic $\mu_1(l_i)$, that is:

$$\lambda_{ext}(l_i \to l_{i+1}) = \left\{ x \in \mu_2(l_i) / \exists t \in \mathbb{R}^+, \mu_1(l_i)(x(t)) \in \lambda(l_i \to l_{i+1}) \right\} \qquad Y$$

We notice that $\lambda_{li \to li+1} \in \lambda ext_{li \to li+1}$.

Explanations about what is X can now be given.

Let consider the path $\pi = (l_1, \ldots, l_i, l_{i+1}, l_{i+2}, \ldots, l_{f-1}, l_f)$. The process initial phase is $l_1$ and a possible process final phase is $l_f$. The phase index i gives the occurrence order of the different phases in the path (and not the phases number of the model). The same phase can occur several times in the path. Now, let define what is the set X:

– if $l_i = l_{f-1}$ then $l_{i+1}$ is the final state $l_f$ of $\pi$ and X is its constraint set: $X = \mu_2(l_f)$

– if $l_i \neq l_{f-1}$ then $l_{i+1}$ is not the final state $l_f$ of $\pi$ and X is the extended jump-set of the two next phases of the path: $X = \lambda_{ext}(l_{i+1} \to l_{i+2})$.

Then comes corollary 1:

*Corollary 1:*

Given the discrete path $\pi$ of the hybrid automata leading from the initial phase to a final phase, given $i \in \{1,\ldots,m\}$, $j \in \{1,\ldots,m\} / l_i \in \pi$, $l_j \in \pi$ are two consecutive discrete phases of $\pi$:

– $\exists i \in \{1,\ldots,m\}$, if $\lambda_{li \to li+1} = \varnothing$, then the constraint $\mu_2$ is not respected along the path and it is not possible to reach the final phase from the initial phase.

– $\forall i \in \{1,\ldots,m\}$, if $\lambda_{li \to li+1} \neq \varnothing$, then the final phase is reachable from the initial phase along the path $\pi$ in a finite time. $\qquad Y$

## 5. Controllability of HDS: application to systems with linear vector fields

### 5.1. *Assumptions and state-space trajectories*

The problem is solved using an analytical approach in the state-space. In the following, except if something else is stated, the system is limited as follow:

*Assumptions:*
− There are only two continuous state variables: $x_1$ and $x_2$
− The system has linear vector fields. Moreover the differential equations described by $\mu_1$ are decoupled, accept converging solutions (which are traditional industrial cases) and are globally Lipschitz vector fields:

$$\mu_1(l_i): \begin{cases} \dot{x}_1 = -a_1.x_1 + b_1 \\ \dot{x}_2 = -a_2.x_2 + b_2 \end{cases} \text{ with } a_j \geq 0$$

− All the real state-variables in all phases must respect a minimum and a maximum value: $\forall j \in \{1,...,n\}, x_{j_{min}} \leq x_j \leq x_{j_{max}}$. Thus the constraint for any phase $\mu_2(l_i)$ is a rectangle defined by horizontal and vertical lines.

First assumption leads to a 2 dimensional state-space problem but the two next theorems are easily extendable to any higher dimension. However it will not be possible to formally characterize these higher dimension systems [ASA 95] that is there are no formal equations which can exactly describe the shape of the regions we will have to deal with.

Third assumption is not a hard limitation because in practical, the constraints for a real variable are often given with a minimum and a maximum value.

In a 2 dimensional state-space, the value of $a_j$, $\forall j \in \{1,2\}$, leads to different cases.

*Theorem 1:*
Given the differential equations described in $\mu_1$:
− if $a_1 = a_2 = 0$, the solutions of the equations are linear and $x_2=f(x_1)$ are straight lines (see figure 4a).
− if $a_1$ or $a_2$ is null but not both ($a_1=0$ or $a_2=0$ exclusively), the set of curves $x_2=f(x_1)$ is composed of exponential curves which converge to a vertical or a horizontal asymptotic line $\Delta$ (see figure 4b). The line equation is: $\Delta: x_i = \dfrac{b_i}{a_i}$
− if $a_1>0$ and $a_2>0$, the set of curves $x_2=f(x_1)$ is composed of exponential curves which converge to an equilibrium point $x_e$ (see figure 4c and 4d), where: $x_{ei} = \dfrac{a_i}{b_i}$ for $i \in \{1, 2\}$

In all cases, the curves $x_2=f(x_1)$ are monotonous.

Y

*Demonstration 1:*

   Let define the differential equation $\dot{x} = -a.x + b$, with $a \geq 0$.
   If $a = 0$, the solution to a such equation is $x(t) = b.t + C$. This is the equation of a line which is monotonous.
   If $a \neq 0$, the solution is $x(t) = \frac{b}{a} + C.e^{-a.t}$. The derivative $\dot{x} = -a.C.e^{-a.t}$. $e^{-a.t} > 0$ and $-a.C = cst$ so the sign of $\dot{x}$ is always the same when t flows. Thus, $x(t) = \frac{b}{a} + C.e^{-a.t}$ is a monotonous curve and if $x_1(t)$ and $x_2(t)$ are monotonous, then $x_2 = f(x_1)$ is also monotonous.
   Therefore, $\forall a \geq 0$, the solution to $\dot{x} = -a.x + b$ is a monotonous curve.

- In case 1, $a_1 = a_2 = 0$ so $x_1(t) = b_1.t + C_1$ and $x_2(t) = b_2.t + C_2$. Thus, $x_2 = \frac{b_2}{b_1}x_1 + \frac{b_2 C_1}{b_1} + C_2$ which is the equation of a straight line.

- In case 2, the differential equation solution is $x(t) = b.t + C$ in a dimension and $x(t) = \frac{b}{a} + C.e^{-a.t}$ in the other dimension. The equation of the curves in the state-space are $x_i = \frac{b_i}{a_i} + C_i.e^{\frac{a_i C_j}{b_j}}.e^{-\frac{a_i}{b_j}x_j}$ with i=1 and j=2 or i=2 and j=1. This is the equation of exponential curves which converge to a horizontal or a vertical asymptotic line $\Delta: x_i = \frac{b_i}{a_i}$ when time flows.

- In case 3, the differential equation solution is $x(t) = \frac{b}{a} + (x_0 - \frac{b}{a}).e^{-a.t}$ in both dimensions. $x_0$ is the value of x when t=0.
   The equation of the curves in the state-space are:

$$x_2 = \frac{b_2}{a_2} + \left(x_{20} - \frac{b_2}{a_2}\right).e^{\frac{a_2}{a_1}\ln\left(\left(x_1 - \frac{b_1}{a_1}\right) \middle/ \left(x_{10} - \frac{b_1}{a_1}\right)\right)}$$

$x_{i0}$ is the initial value of $x_i$.
   This is the equation of exponential curves which converge to the point $(x_1, x_2) = \left(\frac{b_1}{a_1}, \frac{b_2}{a_2}\right)$ when t tends to infinity.

It is reminded that $b_1/a_1$ is the value where $x_1$ converges, $x_{10}$ is the initial value of $x_1$ and the studied curve is monotonous. Therefore, either $x_{10} < x_1 < \dfrac{b_1}{a_1}$ or $x_{10} > x_1 > \dfrac{b_1}{a_1}$, and $\dfrac{x_1 - \dfrac{b_1}{a_1}}{x_{01} - \dfrac{b_1}{a_1}} > 0$. Then, $\ln\left(\dfrac{x_1 - \dfrac{b_1}{a_1}}{x_{01} - \dfrac{b_1}{a_1}}\right)$ is always defined. '

The first case where $a_i=0$ for i=1 and 2 has been studied by Tittus [TIT 93] [TIT 98]. The method described in the following solves cases 2 and 3. In case 3, the relative position of the equilibrium point $x_e$ and the jump-set $\lambda_{li \to li+1}$ leads to two different ways of calculation for the extended jump-set $\lambda ext_{li \to li+1}$.

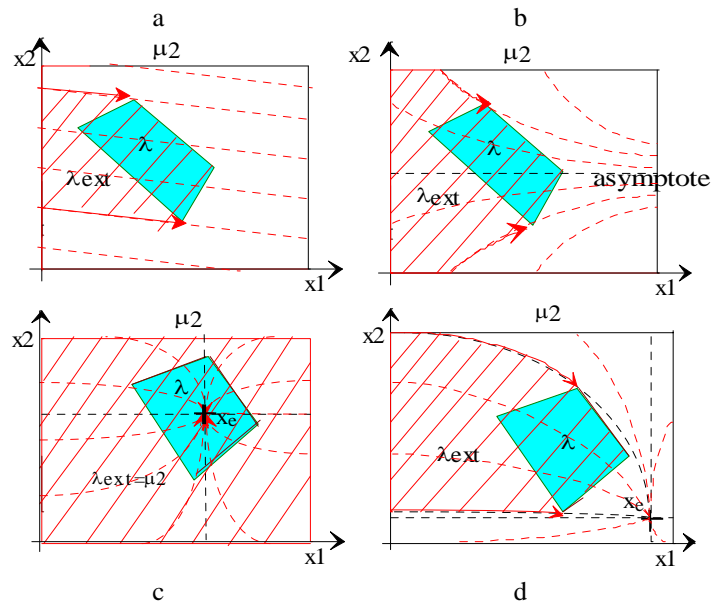

**Figure 4.** *Different cases for λ and λext*

We notice that figure 4 gives a representation of the state-space when there are only two real valued variables ($x \in \mathbb{R}^n$, n=2) but theorem 2 is true for any value of $n \in \mathbb{N}$.

*Theorem 2:*

Given a vector field $\mu_1(l_i)$ of dimension $3^2$ where the equilibrium point $x_e$ is defined, if the equilibrium point $x_e$ is strictly included in the jump-set $\lambda_{li \to li+1}$ then the extended jump-set $\lambda ext_{li \to li+1}$ is bounded by the phase constraint $\mu_2(l_i)$:

$$\lambda ext_{li \to li+1} = \mu_2(l_i)$$

In that case, the equilibrium point can be reached (see figure 4c).     Y

*Demonstration 2:*

In theorem 2, we are in the case where curves converge to an equilibrium point $x_e$. So, it is possible to reach the jump-set from anywhere in the state-space with the dynamic $\mu_1$:

$$\text{if } x_e \in \lambda(l_i \to l_{i+1}), \; \forall x \in \mathbb{R}^2, \exists t \in \mathbb{R}^+ / \; \mu_1(l_i)(x(t)) \in \lambda(l_i \to l_{i+1}) \quad [1]$$

By theorem 1, the curves defined by $\mu_1(l_i)$ are monotonous. Thus, $\dot{x}_1$ and $\dot{x}_2$ are either always positive or always negative.

By assumption, the contour of $\mu_2$ is formed of a set of line segments which are horizontal or vertical. So, $\forall i \in \{1,2\}$, $x_i$=cst and $\dot{x}_i = 0$.

So, the set of curves defined by $\mu_1(l_i)$ will cross at most one time the lines segment defining the contour of $\mu_2(l_i)$.

If the equilibrium point $x_e$ is inside $\mu_2(l_i)$, once a curve defined by $\mu_1(l_i)$ get inside $\mu_2(l_i)$, it will never go out again:

$$\text{if } x_e \in \mu_2(l_i), \; \forall x \in \mu_2(l_i), \forall t \in \mathbb{R}^+, \; \mu_1(l_i)(x(t)) \in \mu_2(l_i) \quad [2]$$

As a consequence of [1] and [2], we get:

$$\text{if } x_e \in \lambda(l_i \to l_{i+1}) \subseteq \mu_2(l_i), \; \forall x \in \mu_2(l_i), \exists t \in \mathbb{R}^+ / \; \mu_1(l_i)(x(t)) \in \lambda(l_i \to l_{i+1})$$

,

The constraints $\mu_2(l_i)$, the jump-set $\lambda_{li \to li+1}$ and the extended jump-set $\lambda ext_{li \to li+1}$ are are the different state-space regions (or sets) we can get. Before we go further, it is necessary to define how these regions can formally be described.

### 5.2. *Analytical expression for a two dimensions region and for its contour*

The aim of this section is to define the existence conditions of an analytical expression which describes the regions we will get.

*Definition 5:*

The contour of the region R is S(R). It is defined by a set of p boundaries $s_k(R)$:

$$S(R)=\{s_k(R), k \in \{1, ..., p\}\}$$

$s_k(R)$ is a curve segment (or boundary) defined by a curve and two extreme points. Each extreme point $ep_j(R)$ belongs to two different boundaries of the contour.of the region. The set of extreme points Pext(R) of a given region R, is defined as follow:

$$P_{ext}(R)=\{ep_j(R), j \in \{1, ..., p\}\}$$

A curve segment $s_k(R)$ is analytically defined by:
- a timed parametered equations system:

$$s_k(R): \begin{cases} x_1 = f_1(t) \\ x_2 = f_2(t) \end{cases}$$

- two extreme points ($ep_a$ and $ep_b$) coordinates:

$$\begin{bmatrix} x_1(ep_a) & x_2(ep_a) \\ x_1(ep_b) & x_2(ep_b) \end{bmatrix} \qquad Y$$

*Definition 6:*

The half-space $Rs_k(R)$ associated to the boundary $s_k(R)$ is the half-space which is over or below the boundary $s_k$. It is analytically defined by:

$$Rs_k(R): \begin{cases} x_1 = f_1(t) \\ x_1(ep_a) \leq x_1 \leq x_1(ep_b) \\ x_2 \geq/\leq f_2(t) \end{cases}$$

$x_2 \geq/\leq f_2(t)$ means that the equation $x_2=f_2(t)$, which is one of the equations defining $s_k(R)$, is transformed into an inequation $x_2 \leq f_2(t)$ or $x_2 \geq f_2(t)$. $\qquad Y$

Now, we can easily imagine that taking the intersection of the half-spaces $Rs_k(R)$ associated to the set of boundaries $s_k(R)$ will define the region R. So an analytical description for these regions seems to be possible. Nevertheless, these regions must verify the property of monotonous regions to be modeled by definition 6 equations / inequations systems.

*Definition 7:*

Let $\sigma$ be a point situated on the contour of a region S(R).

Let $\delta$ be a line which equation is: $x_i$=kY3 where i={1, 2}, k may be equal to 0 but it is not compulsory.

Let define the orthogonal projection from S(R) to $\delta$ which gives an image $\eta \in \delta$ to any point $\sigma$ of S(R):

$$S(R) \xrightarrow{proj} \delta$$
$$\sigma \mapsto \eta$$

Let $proj^{-1}$ be the inverse projection of proj.

$\forall \sigma \in \delta$, if $proj^{-1}(\eta)$ has either an infinite number of images or at most only two images belonging to S(R) then the region R is monotonous. $\qquad Y$

Figure 5 illustrates this property with δ line chosen as a horizontal line. Thus i=2 and its equation is $x_2=k$. For specific points $\eta \in \delta$, the number of image that proj$^{-1}(\eta)$ accepts is counted.

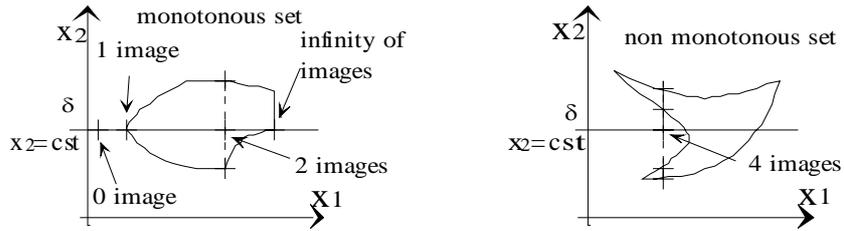

**Figure 5.** *Monotonous set – non monotonous set*

We notice that a monotonous region has not necessarily a contour composed with a set of monotonous boundaries and a set of monotonous boundaries does not always define a monotonous region.

We also notice that a monotonous region can be non convex.

It is time now to realize the link between the boundaries analytical description and the monotonous regions.

*Theorem 3:*

If R is a monotonous region and if its contour is composed of a set of boundaries S(R), then an analytical description for R always exists, being the intersection of the half-spaces $Rs_k(R)$:

$$R = \bigcap_{k=1}^{p} Rs_k(R)$$

and analytically, R is described by:

$$R = \bigcup_{i=1}^{p} \begin{cases} x_{1i} = f_{1i}(t) \\ x_1(ep_{ai}) \leq x_{1i} \leq x_1(ep_{bi}) \\ x_{2i} \geq/\leq f_{2i}(t) \end{cases} \quad Y$$

*Demonstration 3:*

Let R be a monotonous region and $s_k(R) \in S(R)$ be a boundary for R, it always exists an equation system which describes $s_k(R)$. Moreover, it always exists an analytical description of the associated half-space $Rs_k(R)$.

Let $x \in R$ be any point of the monotonous region, x will always belong to two half-spaces associated to S(R) and only two.

So, if we generalize to the set of points x defining the region R, R can have an analytical description which is the union of all the systems describing the half-spaces set:

$$R = \bigcup_{i=1}^{p} \begin{cases} x_{1i} = f_{1i}(t) \\ x_1(ep_{ai}) \leq x_{1i} \leq x_1(ep_{bi}) \\ x_{2i} \geq/\leq f_{2i}(t) \end{cases},$$

Here is an example which gives the analytical description of a monotonous region (figure 6).

$s_1: \begin{cases} x_1(t) = -b_1/a_1 + C_1 \cdot e^{-a_1 \cdot t} \\ x_2(t) \leq -b_2/a_2 + C_2 \cdot e^{-a_2 \cdot t} \\ a \leq x_1 \leq b \end{cases}$

$s_3: \begin{cases} x_1(t) = -b_1/a_1 + C_1 \cdot e^{-a_1 \cdot t} \\ x_2(t) \geq C_2 \\ a \leq x_1 \leq c \end{cases}$

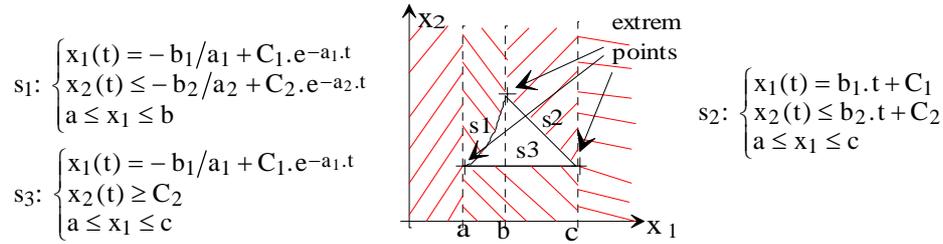

$s_2: \begin{cases} x_1(t) = b_1 \cdot t + C_1 \\ x_2(t) \leq b_2 \cdot t + C_2 \\ a \leq x_1 \leq c \end{cases}$

**Figure 6.** *Monotonous set analytical description*

This analytical method used to describe the regions fits with only 2 dimensional state-spaces. Consequently, the studied systems can have only two state variables.

A n-dimensional monotonous region will have a contour defined by a n-dimensional "hypersurface segments" set. For instance, if we have a 2-dimensional monotonous region, its contour is a 2-dimensional curve segment set, if we have a 3-dimensional monotonous region, its contour is a 3-dimensional hypercurve segment set which would be a kind of curved plane.

But a n-parameter equations system: $x_i = f_i(t)$ where $i = \{1, \ldots, n\}$ describes a single curve in a n-dimensional space. So, it is possible to use this kind of equation systems to describe the "hypersurface segments" only when they are curve segments, i.e. when the space dimension is n=2.

### 5.3. *About trajectory synthesis with two dimensional monotonous regions*

Now we have proven it always exists an analytical description for monotonous regions with contour made of a segment curves set, we want to be sure we will always have to deal with such kinds of regions.

The model we have presented uses three different state-space regions:
- the set of m state constraints $\mu_2(l_i)$
- the jump-set $\lambda_{l_i \to l_{i+1}}$
- the extended-jump-set $\lambda ext_{l_i \to l_{i+1}}$

Firstly, the state constraint of any phase $\mu_2(l_i)$ is rectangular with vertical and horizontal lines by assumption. So it is always monotonous. For the two other sets ($\lambda$ and $\lambda_{ext}$), we must define conditions which are required to be sure they remain monotonous when time flows and when discrete transitions occur.

*Theorem 4:*

The intersection of two monotonous regions is a monotonous region.    Y

*Demonstration 4:*

Let define two monotonous regions $R_a$ and $R_b$ and their contour $S(R_a)$ and $S(R_b)$. It is reminded that a region is monotonous if the orthogonal projection from a horizontal or vertical line to its contour has either an infinite number of images or at most two.

In the current demonstration, the projection line $\delta$ is chosen to be horizontal and never cuts the two regions $R_a$ and $R_b$:

$$\delta = \{\eta \in \delta\, /\, \forall i \in \{a,b\}, \eta \notin S(R_i)\}$$

Let define the two orthogonal projections $proj_a$ and $proj_b$ such that $\forall i \in \{a,b\}$:

$$S(R_i) \xrightarrow{proj_i} \delta$$
$$\sigma_i \mapsto \eta$$

By definition 7, if a region is monotonous then,

$\forall i \in \{a,b\}$, $\forall \sigma_i \in S(R_i)$, $\forall \eta \in \delta\,/\,proj_i^{-1}(\eta) \neq \emptyset$, $proj_i^{-1}(\eta)$ can be:

  – one point $\sigma_{i1}$
  – or two points $\sigma_{i1}$ and $\sigma_{i2}$
  – or an infinite number of points $\sigma_i \in S(R_i)$

If $proj_i^{-1}(\eta)$ is unique then we define a new point $\sigma_{i2} = \sigma_{i1}$.

If $proj_i^{-1}(\eta)$ can be two points, $\sigma_{i1}$ is defined as the furthest point of $\eta$ and $\sigma_{i2}$ as the closest.

If $proj_i^{-1}(\eta)$ is an infinite number of points then we define:

  – $\sigma_{i\,inf} = \{\sigma_i \in S(R_i)\,/\,\sigma_i = proj_i^{-1}(\eta)\}$
  – $\sigma_{i1}$ such that $\sigma_{i1} = \max\limits_{\sigma_i \in \sigma_{i\,inf}}(d(\eta,\sigma_i))$: $\sigma_{i1}$ is the furthest point of $\eta$
  – $\sigma_{i2}$ such that $\sigma_{i2} = \min\limits_{\sigma_i \in \sigma_{i\,inf}}(d(\eta,\sigma_i))$: $\sigma_{i2}$ is the closest point of $\eta$

We notice that $d(x,y)$ means the distance between the two points x and y.

So we can always associate two points ($\sigma_{i1}$ and $\sigma_{i2}$) to every $\eta \in \delta\,/\,proj_i^{-1}(\eta) \neq \emptyset$ for any region $R_i$.

Now, $\forall \eta \in \delta\,/\,proj_a^{-1}(\eta) \neq \emptyset$ and $proj_b^{-1}(\eta) \neq \emptyset$, we must determine which two points of $\sigma_{a1}$, $\sigma_{a2}$, $\sigma_{b1}$ and $\sigma_{b2}$ must be kept in order to sketch the contour of the $R_a$ and $R_b$ intersection. These two points are called $\sigma_1$ and $\sigma_2$.

$\forall \eta \in \delta\,/\,proj_a^{-1}(\eta) \neq \emptyset$ and $proj_b^{-1}(\eta) \neq \emptyset$, we define:

  – $\sigma_1$ such that $d(\delta,\sigma_1) = \min\limits_{k=\{a,b\}}(d(\delta,\sigma_{k1}))$

- $\sigma_2$ such that $d(\delta, \sigma_2) = \max_{k=\{a,b\}}(d(\delta, \sigma_{k2}))$
- $\sigma_{inf}$ such that $\forall i \in \{a, b\}$, if $\sigma_{i\,inf}$ is defined for $\eta$,
  $\sigma_{inf} = \{\sigma_i \in \sigma_{i\,inf} / d(\delta, \sigma_2) \leq d(\delta, \sigma_i) \leq d(\delta, \sigma_1)\}$

If we consider all the possible points $\eta \in \delta / \mathrm{proj}_a^{-1}(\eta) \neq \varnothing$ and $\mathrm{proj}_b^{-1}(\eta) \neq \varnothing$, a set of points $\sigma_1$, $\sigma_2$ and $\sigma_{inf}$ is defined. This set will define the contour $S(R)$ of the region $R$ / $R = R_a \cap R_b$.

Finally, let define $S(R) \xrightarrow{\mathrm{proj}} \delta$ and the inverse projection $\mathrm{proj}^{-1}$, $\forall \eta \in \delta$, $\mathrm{proj}_i^{-1}(\eta)$ will be:
- not defined
- or one point: $\sigma_1 = \sigma_2$
- or two points: $\sigma_1 \neq \sigma_2$
- or an infinite number of points $\in \sigma_{inf}$

Thus $S(R)$ has either an infinite number of images or at most two images. This is the definition of a monotonous region.

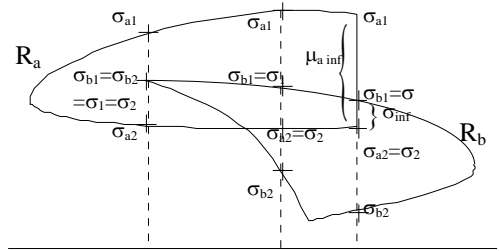

**Figure 7.** *Intersection of two monotonous regions*

*Theorem 5:*

For a given couple of location $l_i$ and $l_{i+1}$, if $\mu_2(l_i)$ and $\lambda \mathrm{ext}_{l_{i+1} \to X}$ are monotonous regions and if $\mu_3(l_i, \delta_{i,i+1}, l_{i+1})$ conserves the monotonous property, then $\lambda_{l_i \to l_{i+1}}$ is monotonous.

*Demonstration 5:*

If $\lambda \mathrm{ext}_{l_{i+1} \to X}$ is monotonous and if $\mu_3(l_i, \delta_{i,i+1}, l_{i+1})$ conserves the monotonous property then the set $\nu$ such that if $x \in \nu$ then $\mu_3(l_i, \delta_{i,i+1}, l_{i+1})(x) \in \lambda \mathrm{ext}_{l_{i+1} \to X}$ is a monotonous region.

Definition 3 says that $\lambda_{l_i \to l_{i+1}} = \nu \cap \mu_2(l_i)$. By assumption, $\mu_2(l_i)$ is monotonous and by theorem 4, $\lambda_{l_i \to l_{i+1}}$ is monotonous too.

The monotonous property conservation for function $\mu_3$ is not well defined but we can note that any action on a state–variable must have an equivalent effect on all the others.

For instance, let consider a two dimensional state-space, i.e. $[x_1, x_2]^T = x \in \mathbb{R}^2$. The function $\mu_3$ such that $x_1(t^+)=x_1(t^-)+k_1$ and $x_2(t^+)=x_2(t^-)+k_2$ is an example of what is an equivalent effect for both variables.

Moreover, we are sure that if there are no continuous discontinuity when a transition occurs, that is if $\mu_3(l_i \to l_{i+1})=\text{Id}$, then $\mu_3$ conserves the monotonous property. This is the case in the example presented in section 2 and solved in section 6.

*Definition 8:*

If $\lambda ext_{l_i \to l_{i+1}} \neq \mu_2(l_i)$, the extended jump-set contour ($\lambda ext_{l_i \to l_{i+1}}$ contour) is composed of three kinds of boundaries (or curve segments):
- some $\lambda_{l_i \to l_{i+1}}$ boundaries
- some $\mu_2(l_i)$ boundaries
- curve segments which are solutions of the differential equations described by $\mu_1(l_i)$. These curves are called extension curves. they can cross a $\lambda$-boundary or be a tangente to $\lambda$                                                                           Y

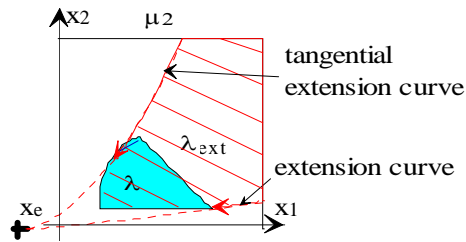

**Figure 8.** *Extension curves and $\lambda ext$-contour*

Now we have seen how the contour of the extended jump-set $S(\lambda ext)$ is composed, we can easily understand that for a given couple of location $l_i$ and $l_{i+1}$, if $\lambda_{l_i \to l_{i+1}}$ is monotonous and if the state-variables are not coupled (i.e. the solution to $\dot{x} = A.x + b$ is a monotonous curve), then $\lambda ext_{l_i \to l_{i+1}}$ is monotonous.     Y

*Corollary 2:*

If the set of m state constraints $\mu_2$ of any phase $l_i$ is always monotonous, if $\mu_3(l_i, \delta_{i,i+1}, l_{i+1})$ conserves the monotonous property and if the evolution function $\mu_1(l_i)$ is such that the state-variables are not coupled, then all the regions we will get when computing the reachability set will be monotonous. So, a formal description of all encountered sets will always be possible.                                                            Y

Corollary 2 defines the conditions we need to be able to formally compute the reachable set using an algorithm.

## 6. Algorithm and application to the batch system

### 6.1. *λext computation algorithm for two dimensional state-space problems*

The hardest step to decide whether a system is controllable or not is the computation of the extended jump-set λext for a given location $l_i$. The algorithm is given below:

- *initialization of:*
    - *the evolution function $\mu_1$ modeled by differential equations*
    - *the jump-set $\lambda$ and the constraint set $\mu_2$ by their analytical description*
- *if the equilibrium point is defined and is strictly included in the jump-set $\lambda$ then:*
    $\lambda ext = \mu_2$
  *else:*
    - *calculation of all extreme-points*
    - *calculation of all the curves having $\mu_1(l_i)$ dynamic trajectory and crossing a $\lambda$ extreme-point*
    - *calculation of $\lambda$ tangential trajectory curves if they exist*
    - *calculation of the extension curves by selecting the curves defining the largest region*
    - *appropriate $\lambda$ and $\mu_2(l_i)$ boundaries are added to the extension curves to form $\lambda_{ext}(l_i \rightarrow l_{i+1})$ contour*
    - *equalities are transformed into inequalities to analytically define $\lambda_{ext}(l_i \rightarrow l_{i+1})$*

More details are given about this algorithm in [MAN99a]. We remind that in the current paper, the goal is to present the analytical definition of regions which are necessary to the calculation of the extended jump-set for a given phase.

### 6.2. *Application to the batch system*

Let us apply this control synthesis methodology to the batch system example described in section 2 (see figure 1) and modeled in section 3 (see figure 2).

A discrete controllable path π must be calculated. It must start from the initial state $l_0$ and stop at the final state $l_3$. The path computation is realized backward. The jump-sets, extended jump-sets and phase constraints are represented in the state-space $x_2=f(x_1)$. A new state-space representation is given for each phase the system goes through.

A controllable discrete path π is synthesized:

$$\pi = (l_0, l_1, l_2, l_1, l_3)$$

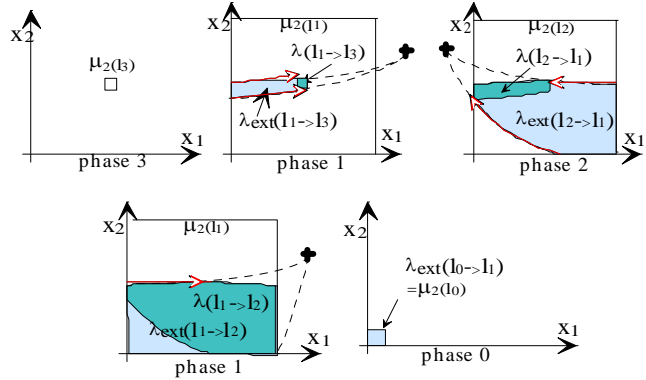

**Figure 9.** *Controllability analysis realized in the state-space*

## 7. Concluding remarks

We have presented parts, dealing with region analytical description, of a new method which goal is to synthesize a return trajectory from a default mode to one of the nominal modes for hybrid dynamical systems with linear vector fields in a two dimensional state-space. It is an extension of Tittus's works on integrator hybrid systems. A great care has been taken in proving the computability of the state-space regions involved into the reachability analysis.

Firstly, we have presented a brief example with a batch system in a default mode to be driven back to its nominal mode.

Once the position of the problem was clearly fixed, some theoretical reminds from Tittus's researches, necessary to solve the problem, were presented.

Afterward, we have explained it is always possible to formally describe a monotonous region and we have defined the necessary conditions to always have such kind of regions when computing a safe hybrid trajectory.

In the last section, the algorithm which computes the extended jump-set for a given location was given and a solution for the batch system example was presented.

At that point of research, the assumptions made for the continuous state variables $x_i$ are rather strong and limit the possible applications field. One next step is to relax them. Particularly systems with coupled variables are now under investigation and first results are positive [CHA 99].

Once the jump-sets and extended jump-sets are computed for a given path, it is necessary to choose a control strategy (e.g. switch as soon as possible, switch as late as possible, …) in order to determine the switching instants.

An another field under investigation is to develop a methodology which can achieve an optimal control strategy. If the criteria to be optimized is time independent, then the control sequence design seems to be computable using pure discrete optimization tools. In the other case, continuous and discrete optimization tools must be combined.